\documentclass[namedreferences]{kluwer}

   \def\kaplogo{}\def\thecopyright{}
   \def\BibTeX{{\rmfamily B\kern-.05em \textsc{i\kern-.025em b}%
   \kern-.08em T\kern-.1667em\lower.7ex\hbox{E}\kern-.125emX}}

\usepackage[dvips]{graphics}

\newtheorem{theorem}{Theorem}[section]
\newtheorem{definition}[theorem]{\it Definition}

\newcommand{\Proof}[1][{\bf Proof.}]{\begin{proof}[#1]}
\newcommand{\QED}{\end{proof}}

\newcommand{\EQn}[1]{{\rm (\ref{EQ:#1})}}

\renewcommand{\theequation}{\thesection.\arabic{equation}}

\newcommand{\Ac}{\mathcal{A}}
\newcommand{\Cc}{\mathcal{C}}
\newcommand{\Ec}{\mathcal{E}}
\newcommand{\Fc}{\mathcal{F}}
\newcommand{\Ic}{\mathcal{I}}
\newcommand{\Lc}{\mathcal{L}}
\newcommand{\Mc}{\mathcal{M}}
\newcommand{\Sc}{\mathcal{S}}

\newcommand{\PP}{{\bf P}}
\newcommand{\UU}{{\bf U}}
\newcommand{\XX}{{\bf X}}

\newcommand{\Prob}{\mathbb{P}}

\newcommand{\stocle}{\preceq}

\newcommand{\asection}[1]{\lceil{#1}\rceil}

\begin{document}

\begin{opening}
   \title{Realizable Monotonicity and Inverse Probability Transform}
   \author{
       James Allen Fill\email{jimfill@jhu.edu}
       \thanks{Research for both authors was supported by
          NSF grant DMS-98-03780.}}
   \institute{Department of Mathematical Sciences,
     The Johns Hopkins University}
   \author{
       Motoya Machida\email{machida@math.usu.edu}}
   \institute{Department of Mathematics and Statistics,
     Utah State University}

   \date{July 19, 2000}

   \begin{abstract}
A system $(P_\alpha: \alpha \in A)$ of probability measures on a
common state space $S$ indexed by another index set $A$
can be ``realized'' by a system $(X_\alpha: \alpha \in A)$ of $S$-valued
random variables on some probability space
in such a way that each $X_\alpha$ is distributed as $P_\alpha$.
Assuming that $A$ and $S$ are both partially ordered,
we may ask when the system $(P_\alpha: \alpha \in A)$
can be realized by a system $(X_\alpha: \alpha \in A)$ with
the monotonicity property
that $X_\alpha \le X_\beta$ almost surely whenever $\alpha \le \beta$.
When such a realization is possible, we call the system
$(P_\alpha: \alpha \in A)$ ``realizably monotone.''
Such a system necessarily is stochastically monotone,
that is, satisfies $P_\alpha \le P_\beta$ in stochastic ordering whenever $\alpha \le \beta$.
In general, stochastic monotonicity is
not sufficient for realizable monotonicity.
However, for some particular choices of partial orderings in a finite state setting,
these two notions of monotonicity are equivalent.
We develop an inverse probability transform
for a certain broad class of posets $S$,
and use it to explicitly construct
a system $(X_\alpha: \alpha \in A)$ realizing the monotonicity
of a stochastically monotone system when the two notions of monotonicity are equivalent.
   \end{abstract}


   \classification{AMS subject classification}{Primary 60E05; 
   secondary 06A06, 60J10, 05C05, 05C38.}

\end{opening}



\section{Introduction}

\subsection{Two notions of monotonicity}

We will discuss two notions of monotonicity for probability measures on 
a finite partially ordered set (poset).
Let $\Sc$ be a finite poset and
let $(P_1,P_2)$ be a pair of probability measures on $S$.
(We use a calligraphic letter $\Sc$ in order to distinguish the set
$S$ from the same set equipped with a partial ordering $\le$.)
A subset $U$ of $S$ is said to be an {\em up-set\/} in $\Sc$
(or {\em increasing set\/})
if $y \in U$ whenever $x \in U$ and $x \le y$.
We say that $P_1$ is {\em stochastically smaller\/} than $P_2$,
denoted $P_1 \stocle P_2$,
if
\begin{equation}\label{EQ:stocle.def}
P_1(U) \le P_2(U)
\quad\mbox{ for every up-set $U$ in $\Sc$.}
\end{equation}
An important characterization of stochastic ordering was established
by~\inlinecite{Strassen} and fully investigated by~\inlinecite{KKO}.
They show that \EQn{stocle.def} is necessary and sufficient
for the existence of a pair $(\XX_1, \XX_2)$ of $S$-valued random variables
[defined on some probability space $(\Omega,\Fc,\Prob)$] satisfying
the properties that $\XX_1 \le \XX_2$ and that
\linebreak
$\Prob(\XX_i \in \cdot) = P_i(\cdot)$ for $i=1,2$.

Now let $\Ac$ be a finite poset.
Let $(P_\alpha: \alpha\in A)$ be a system of probability measures on $S$.
We call $(P_\alpha:\alpha\in{A})$ a {\em realizably monotone\/} system
if there exists a system $(\XX_\alpha:\alpha\in A)$
of $S$-valued random variables
such that
\begin{equation}\label{EQ:rm.mono}
  \XX_\alpha \le \XX_\beta
  \quad\mbox{ whenever $\alpha \le \beta$ }
\end{equation}
and
\begin{equation}\label{EQ:rm.marg}
\Prob(\XX_\alpha \in \cdot) = P_\alpha(\cdot)
\quad\mbox{ for every $\alpha\in A$. }
\end{equation}
In such a case we shall say that $(\XX_\alpha:\alpha\in{A})$
{\em realizes the monotonicity\/} of $(P_\alpha:\alpha\in{A})$.
The (easier half of the) characterization of stochastic ordering applied pairwise implies
\begin{equation}\label{EQ:sm.mono}
  P_\alpha \stocle P_\beta \quad\mbox{ whenever $\alpha\le\beta$. }
\end{equation}
The system $(P_\alpha: \alpha\in A)$ is said to be {\em stochastically monotone\/}
if it satisfies~\EQn{sm.mono}.
Thus, stochastic monotonicity is necessary for realizable
monotonicity.

In light of Strassen's characterization of stochastic ordering,
one might guess that stochastic monotonicity is also sufficient for
realizable monotonicity.
It is perhaps surprising that the conjecture is false in general.
Various counterexamples are given by~\inlinecite{SMRM},
including one independently discovered by~\inlinecite{Ross}.
Given a pair $(\Ac,\Sc)$ of posets,
if the two notions of monotonicity---stochastic and realizable---are
equivalent, then we say that {\em monotonicity equivalence\/} holds
for $(\Ac,\Sc)$.

\subsection{Inverse probability transform}\label{intro:inv}

Suppose that $\Sc$ is linearly ordered.
Then, for a given probability measure $P$ on $S$,
we can define its {\em inverse probability transform\/}
$P^{-1}$ by
\begin{equation}\label{inv}
  P^{-1}(t)
  := \min\left\{x\in S: t < F(x) \right\},
  \quad t \in [0,1),
\end{equation}
where $F(x) := P(\{z \in S: z \le x\})$ is
the distribution function for~$P$.
Furthermore, let $\Ac$ be any poset, and
let $(P_\alpha: \alpha\in A)$ be a stochastically monotone system of
probability measures on $S$.
Given a single uniform random variable $\UU$ on $[0,1)$,
we can construct a system $(\XX_\alpha: \alpha\in A)$ of
$S$-valued random variables via
$\XX_\alpha := P_\alpha^{-1}(\UU)$
which realizes the monotonicity.
This proves that
monotonicity equivalence always holds for $(\Ac,\Sc)$
when $\Sc$ is linearly ordered.

In Section~\ref{SS:invprob} we generalize the definition of inverse probability transform
to a certain class of posets~$\Sc$ which are not necessarily linearly
ordered.
We then extend the construction in the preceding paragraph
and present Theorems~\ref{TH:class.z} and~\ref{TH:class.w},
thereby establishing monotonicity
equivalence under certain additional assumptions.
A further extension of Theorem~\ref{TH:class.w} is discussed briefly
in Section~\ref{SS:me.class.w}, which
culminates in Theorem~\ref{TH:class.w.ext}.
We will not discuss the proofs of Theorems~\ref{TH:class.z},
\ref{TH:class.w}, and~\ref{TH:class.w.ext} in the present brief paper,
but rather refer the reader to~\inlinecite{thesis}
for (the highly technical) proofs and more extensive discussion.

\subsection{Importance in perfect sampling algorithms}

Of particular interest in our general study of realizable monotonicity
is the case $\Ac = \Sc$.
Here the system $(P(x,\cdot): x \in S)$
of probability measures can be considered as a Markov
transition matrix $\PP$ on the state space $S$.
\inlinecite{Propp-Wilson} and \inlinecite{Fill}
introduced algorithms to produce observations distributed
{\em perfectly\/} according to the long-run distribution of a Markov
chain.
Both algorithms apply most readily and operate most efficiently when
the state space $\Sc$ is a poset and a suitable monotonicity condition
holds.
Of the many differences between the two algorithms, one is that the
appropriate notion of monotonicity for the Propp--Wilson algorithm is
realizable monotonicity, while for Fill's algorithm it is stochastic
monotonicity; see Remark~4.5 in~\inlinecite{Fill}.
Here the properties~\EQn{rm.mono}--\EQn{rm.marg} are essential
for the Propp--Wilson algorithm to be able to generate
transitions simultaneously from every state in such a way
as to preserve ordering relations.
For further discussion of these perfect sampling algorithms in the monotone setting,
see~\inlinecite{Fill} and~\inlinecite{Propp-Wilson};
for further discussion of perfect sampling in general, consult the annotated
bibliography at {\tt http://www.dbwilson.com/exact/}.
\inlinecite{SMRM} show that
the two notions of monotonicity are equivalent if and only if the poset $\Sc$ is acyclic;
see Section~\ref{SS:poset} herein for the definition of this term.

\section{A generalization of inverse probability transform}\label{SS:invprob}

\subsection{Distribution functions on an acyclic poset}\label{SS:poset}

We begin with a notion of acyclic poset, and
its use in introducing a distribution function
on such a poset.
Most of the basic poset terminology adopted here can be found
in~\inlinecite{Stanley} or~\inlinecite{Trotter},
and most of the graph-theoretic terminology in~\inlinecite{West}.
Let $\Sc$ be a poset.
For $x, y \in S$,
we say that $y$ {\em covers\/} $x$ if $x < y$ in $\Sc$ and no element
$z$ of $S$ satisfies $x < z < y$.
We define the {\em cover graph\/} $(S,\Ec_{\Sc})$ of $\Sc$
to be the undirected graph with edge set $\Ec_{\Sc}$
consisting of those unordered pairs $\{x,y\}$
such that either $x$ covers $y$ or $y$ covers $x$ in $\Sc$.
A poset $\Sc$ is said to be {\em acyclic\/} when its cover
graph $(S,\Ec_{\Sc})$ is acyclic in the usual graph-theoretic sense
(i.e., the graph has no cycle).

Throughout the sequel we assume that
the cover graph $(S,\Ec_{\Sc})$ is acyclic and also connected,
that is, that the graph $(S,\Ec_{\Sc})$ is a tree.
Let $\tau$ be a fixed leaf of $(S,\Ec_{\Sc})$,
that is, an element $\tau$ in $S$ such that there exists
a unique edge $\{\tau,z\}$
in $\Ec_{\Sc}$ (for some $z \in S$).
Then, declare $x \le_{\tau} y$ for $x,y \in S$
if the (necessarily existent and unique)
path $(\tau,\ldots,x)$ in the graph from $\tau$ to $x$ contains
the path $(\tau,\ldots,y)$ from $\tau$ to $y$ as a segment.
This introduces a partial ordering $\le_{\tau}$ on the ground set
$S$ [\inlinecite{Bogart}], which may be different from~$\leq$ for the original poset $\Sc$.
We call this new poset $(S,\le_{\tau})$ a {\em rooted tree\/} (rooted
at $\tau$).

For each $x \in S$, set
$$
  \Cc(x) := \{z \in S: \mbox{$x$ covers $z$ in $(S,\le_{\tau})$}\}.
$$
Then a linear extension $(S,\le_{\psi})$ of $(S,\le_{\tau})$
can be obtained by choosing
a linear ordering on $\Cc(x)$ for every $x \in S$.
Explicitly, we define $x \le_{\psi} y$ if either
(i) $x \le_{\tau} y$, or
(ii) there exist some $z \in S$ and some $w, w' \in \Cc(z)$
such that $x \le_{\tau} w$, $y \le_{\tau} w'$,
and $w$ has been chosen to be smaller than $w'$ in $\Cc(z)$.
See Section~\ref{SS:example} for an example.

\begin{definition}{\em
For a given probability measure $P$ on $S$,
we define the {\em distribution function\/} $F(\cdot)$ of $P$ by
$$
  F(x) := P(\{z \in S: z \le_{\tau} x\})
  \quad\mbox{ for each $x \in S$,}
$$
and the {\em distribution function $F\asection{\cdot}$ of linear extension\/} by
$$
  F\asection{x} := P(\{z\in S: z \le_{\psi} x\})
  \quad\mbox{ for each $x \in S$.}
$$
}\end{definition}

In particular, when $(S,\Ec_{\Sc})$ is a path from one end point
$\tau$ to the other end point,
the rooted tree $(S,\le_{\tau})$ is linearly ordered,
and therefore $F(\cdot) \equiv F\asection{\cdot}$.

\subsection{Inverse probability transform}

For a given distribution function $F\asection{\cdot}$ of linear extension
on $S$, we define the
{\em inverse probability transform\/} $P^{-1}$,
a map from $[0,1)$ to $S$, by
\begin{equation}\label{g.inv}
P^{-1}(t) :=
\min\{x \in S: t < F\asection{x}\,\}
\quad\mbox{ for $t \in [0,1)$, }
\end{equation}
where the minimum is given in terms of the linearly ordered set
$(S,\le_{\psi})$.
When $\Sc$ is linearly ordered,
the two definitions of inverse probability transform
in~(\ref{inv}) and~(\ref{g.inv}) are the same.
This equivalence can be extended to the case
that the cover graph $(S,\Ec_{\Sc})$ is a path,
because of the fact that then $F(\cdot) \equiv F\asection{\cdot}$.
Moreover, the property of inverse probability transform
discussed in Section~\ref{intro:inv} remains true in that case:
\begin{theorem}\label{TH:class.z}
Let $\UU$ be a 
uniform random variable on $[0,1)$.
Suppose that $(S,\Ec_{\Sc})$ is a path.
Then, a stochastically monotone system
$(P_\alpha: \alpha\in A)$ is always realizably monotone
via $\XX_\alpha := P_\alpha^{-1}(\UU)$.
\end{theorem}
Theorem~\ref{TH:class.z} reiterates a result presented by
\inlinecite{SMRM}, namely, Theorem~6.1 in their paper.
An acyclic poset $\Sc$ is called a {\em poset of Class\/}~Z
if the cover graph $(S,\Ec_{\Sc})$ is a path.
Otherwise, the acyclic (connected) poset $\Sc$ has a multiple-element $\Cc(x)$
for some $x \in S$.
An example in Section~\ref{SS:example} will demonstrate that
Theorem~\ref{TH:class.z} can fail when an acyclic poset $\Sc$ is not in
Class~Z.
Besides the result for Class~Z, Fill and Machida
gave a complete answer to the monotonicity equivalence problem
[i.e., the question whether monotonicity equivalence holds for given $(\Ac, \Sc)$]
when there exists some $x \in S$ such that (i) $\Cc(x)$ contains
at least two elements, and (ii) $x$ is neither minimal nor maximal in
$\Sc$, that is, when an acyclic poset $\Sc$ falls into
either {\em Class\/}~B or {\em Class\/}~Y, in their terms.
[In their investigation,
a construction of random variables with the desired
properties~(\ref{EQ:rm.mono})--(\ref{EQ:rm.marg}) was reduced to
application of Strassen's characterization of stochastic ordering
if monotonicity equivalence holds
for $(\Ac,\Sc)$ with $\Sc$ a poset either of Class~B or of Class~Y.]
However, when $\Sc$ is a poset satisfying the property that
$x$ is either maximal or minimal in $\Sc$ whenever
$\Cc(x)$ contains at least two elements,
which they (and we) call a {\em poset of Class\/}~W,
we do not know a complete answer
to the monotonicity equivalence problem.
But for a poset~$\Sc$ of Class~W
our generalization of inverse probability transform can,
for some posets~$\Ac$, be used to establish monotonicity equivalence:
\begin{theorem}\label{TH:class.w}
Let $\UU$ be a
uniform random variable on $[0,1)$.
Suppose that $\Sc$ is a poset of Class~{\rm W}, and that
$\Ac$ is a poset having a minimum element and a maximum element.
Then, given a stochastically monotone system $(P_\alpha: \alpha\in A)$,
there exists a system $(\phi_\alpha: \alpha\in A)$
of $\UU$-invariant maps from $[0,1)$ to $[0,1)$ [i.e., $\phi(\UU) \stackrel{\Lc}{=} \UU$]
such that
$$
  \XX_\alpha := P_\alpha^{-1}(\phi_\alpha(\UU)),
  \quad \alpha\in A,
$$
realizes the monotonicity.
\end{theorem}
We call the $\UU$-invariant maps $\phi_\alpha$
in Theorem~\ref{TH:class.w}
{\em synchronizing functions\/}.
Discussion about
how we can practically construct the desired synchronizing functions
can be found in \inlinecite{thesis}.

\subsection{An example}\label{SS:example}

Consider the poset $\Sc$ of Class~W with the following Hasse diagram:
$$
\Sc :=
\begin{minipage}{1.6in}\begin{center}
  \includegraphics{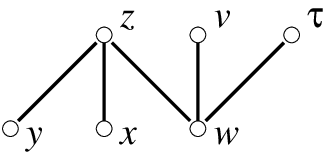}
\end{center}\end{minipage};
$$
the downward arc from $z$ to $y$, for example,
indicates that $z$ covers $y$ in $\Sc$.
Then the rooted tree $(S,\le_{\tau})$ rooted at $\tau$
and the corresponding linear extension $(S,\le_{\psi})$
are given respectively by
$$
(S,\le_{\tau}) =
\begin{minipage}{1.4in}\begin{center}
  \includegraphics{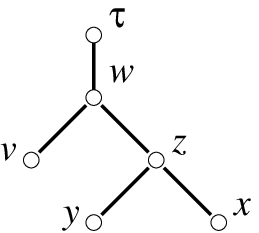}
\end{center}\end{minipage}
\mbox{ and }\quad
(S,\le_{\psi}) =
\begin{minipage}{0.8in}\begin{center}
  \includegraphics{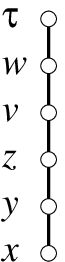}
\end{center}\end{minipage},
$$
where we have chosen the linear ordering $z < v$ in $\Cc(w)$
and the linear ordering $x < y$ in $\Cc(z)$.
Let $P_1$ and $P_2$ be the two probability measures on $S$ in the following table:
$$
\begin{array}{l|cccccc}
\xi & x & y & z & v & w & \tau \\ \hline
P_1(\{\xi\}) & 3/15 & 2/15 & 1/15 & 1/15 & 7/15 & 1/15 \\
P_2(\{\xi\}) & 1/15 & 1/15 & 6/15 & 3/15 & 2/15 & 2/15
\end{array}
\quad.
$$

\begin{figure}[h]
\includegraphics{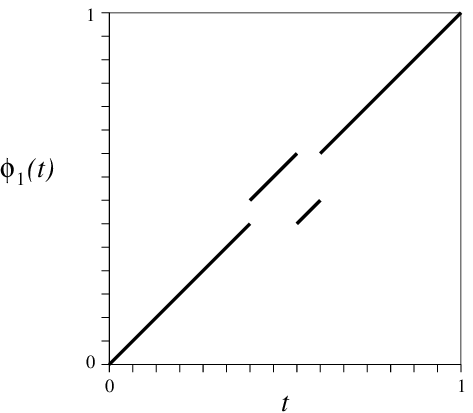}
\hspace{0.2in}
\includegraphics{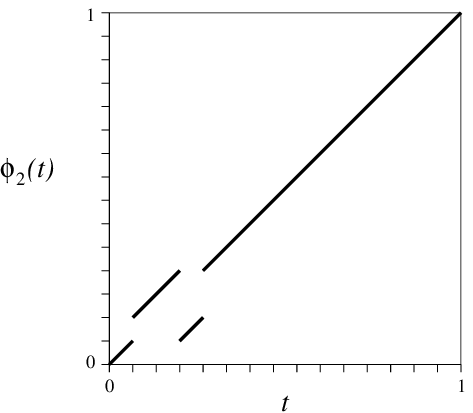}

\caption{The synchronizing functions $\phi_1$ and $\phi_2$}
\label{phi12.fig}
\end{figure}

We can easily check that $P_1 \stocle P_2$.
However,
$P_1^{-1}(t) \le_{\Sc} P_2^{-1}(t)$ does not hold for all $t \in [0,1)$:
For $t \in [\frac{1}{15},\frac{2}{15})$,
$x = P_1^{-1}(t)$ is incomparable with $y = P_2^{-1}(t)$,
and again for $t \in [\frac{6}{15},\frac{7}{15})$,
$v = P_1^{-1}(t)$ is incomparable with $z = P_2^{-1}(t)$.

\begin{figure}[h]
\includegraphics{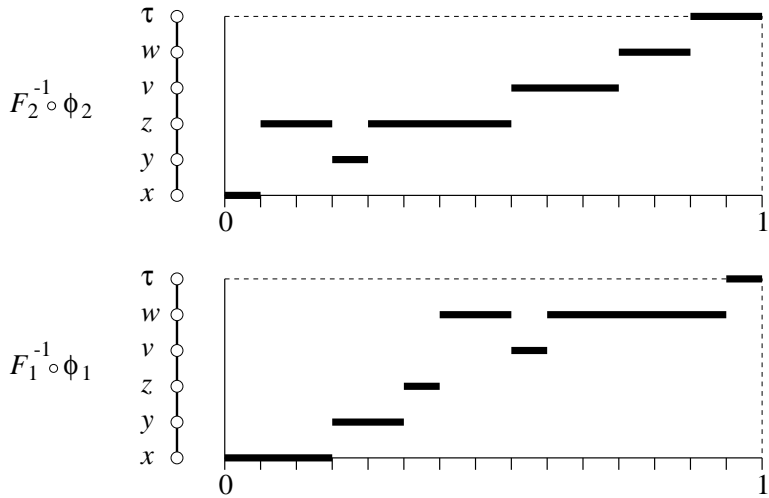}

\caption{The synchronized inverse probability transforms
$P_1^{-1}\circ\phi_1$ and $P_2^{-1}\circ\phi_2$}
\label{Fphi12c.fig}
\end{figure}

Figure~\ref{phi12.fig} displays the synchronizing functions $\phi_1$ and $\phi_2$,
which are both $\UU$-invariant from $[0,1)$ to $[0,1)$.
Then consider the map $P_k^{-1}\circ\phi_k$ from $[0,1)$ to $S$
for $k=1$ and $k=2$, as in Figure~\ref{Fphi12c.fig}.
It is clear from Figure~\ref{Fphi12c.fig} that
$P_1^{-1}(\phi_1(t)) \le_{\Sc} P_2^{-1}(\phi_2(t))$
for all $t \in [0,1)$, as desired.

\section{More on the monotonicity problem}\label{SS:me.class.w}

\subsection{Synchronizable posets}

Given a poset $\Sc$ of Class~W,
Theorem~\ref{TH:class.w} implies
that if $\Ac$ is a poset having
a minimum element and a maximum element,
then monotonicity equivalence holds for $(\Ac,\Sc)$.
This section introduces without detail
a further extension of monotonicity equivalence to a synchronizable
poset~$\Ac$ (as defined below).

Let $D_{\Ac}$ be the set of all the minimal elements in $\Ac$.
Then we define a graph $(D_{\Ac},\Ic_{\Ac})$ on the vertex set $D_{\Ac}$
by including $\{\alpha,\alpha'\}$ as an edge in $\Ic_{\Ac}$
if $\alpha \neq \alpha'$ and
there exists some $\beta \in A$ such that $\alpha,\alpha' < \beta$
in $\Ac$.
We define in analogous fashion a graph $(D^*_{\Ac},\Ic^*_{\Ac})$ on
the set $D^*_{\Ac}$ of all the maximal elements in $\Ac$.
We call these {\em graphs of interlacing relation\/}.
Let $(D_{\Ac},\Ic_0)$ be a spanning tree of $(D_{\Ac},\Ic_{\Ac})$,
that is, let $(D_{\Ac},\Ic_0)$ be a tree with $\Ic_0 \subseteq \Ic_{\Ac}$.
We will say that
$(D_{\Ac},\Ic_0)$ is a {\em locally connected spanning tree\/} of
$(D_{\Ac},\Ic_{\Ac})$ if for every $\alpha \in A$,
the subgraph of $(D_{\Ac},\Ic_0)$ induced by
$D_{\Ac}(\alpha)
:= \{\beta \in D_{\Ac}:
   \mbox{ $\beta\le\alpha$ in $\Ac$}\}$
is connected.
Finally, we call $\Ac$ a {\em synchronizable poset\/}
if there exist respective locally connected spanning trees
of $(D_{\Ac},\Ic_{\Ac})$ and of $(D^*_{\Ac},\Ic^*_{\Ac})$.

\begin{theorem}\label{TH:class.w.ext}
If $\Sc$ is a poset of Class~{\rm W} and $\Ac$ is a synchronizable poset,
then monotonicity equivalence holds for $(\Ac,\Sc)$.
\end{theorem}

Theorem~\ref{TH:class.w.ext} [which is Theorem~6.2 in~\inlinecite{thesis}]
is the most general positive result we know for the monotonicity equivalence
problem when~$\Sc$ is a poset of Class~W.

\subsection{Open problem}

Let
$$
\Sc :=
\begin{minipage}{1.2in}\begin{center}
  \includegraphics{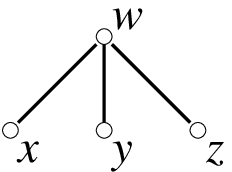}
\end{center}\end{minipage},
$$
which is a poset of Class~W, and let
$$
\Ac :=
\begin{minipage}{1.8in}\begin{center}
  \includegraphics{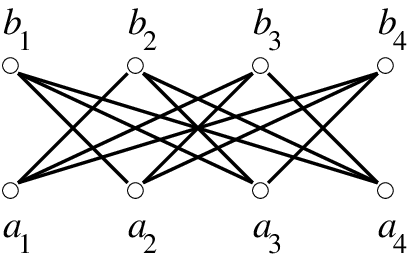}
\end{center}\end{minipage}.
$$
Then the poset $\Ac$ is not synchronizable.
However, we can show that monotonicity equivalence holds for
$(\Ac,\Sc)$ [\inlinecite{thesis}].

Theorem~\ref{TH:class.w.ext} has shown that 
synchronizability of the poset $\Ac$ is sufficient
for monotonicity equivalence when $\Sc$ is a poset of Class~W.
But the above example disproves the assertion that synchronizability
is necessary for monotonicity equivalence.
Furthermore, let $\Mc(\Sc)$ denote the class of all posets $\Ac$ of
monotonicity equivalence for $\Sc$.
Then we can also demonstrate
[cf.\ Example~6.33 in~\inlinecite{thesis}]
that $\Mc(\Sc)$ is {\em not\/} the same
for all posets $\Sc$ of Class~W.
Thus, the interesting question raised but not settled by the present paper
is how to completely characterize posets $\Ac$ of monotonicity
equivalence given a poset $\Sc$ of Class~W,
that is, to determine $\Mc(\Sc)$
exactly for each poset $\Sc$ of Class~W.

\vspace{0.2in}

{\bf Acknowledgments.}
The second author carried out research leading to this paper while he
was a doctoral student in the Department of Mathematical Sciences at
the Johns Hopkins University.
We thank Keith Crank, Alan Goldman, Leslie Hall,
and Edward Scheinerman for providing helpful comments.

%
%
%
%
%
%
%
%
%
%
%
%
%
%
%


\end{document}